\documentclass[10pt]{article}
\usepackage{amsmath,amssymb,amsthm,amscd}
\numberwithin{equation}{section}

\def\RR{{\mathbb R}}

\newtheorem{prop}{Proposition}[section]
\newtheorem{theo}[prop]{Theorem}
\newtheorem{lemm}[prop]{Lemma}
\newtheorem{coro}[prop]{Corollary}
\newtheorem{rema}[prop]{Remark}

\def\begeq{\begin{equation}}
\def\endeq{\end{equation}}

\def\lab{\ }
\def\lab{\label}

\begin{document}
%\tableofcontents
\title{Rigidity of Asymptotically Hyperbolic Manifolds}
%\thanks{}
\author{Yuguang Shi \thanks{The first author's research  is  partially supported by NSF of China,
project 10231010.}\\
Key Laboratory of Pure and Applied Mathematics
\\School of Mathematics Science,\\
    Peking University\\
 Beijing, 100871, China \\
 email:\texttt{ygshi@math.pku.edu.cn}\\[3pt]
 \and Gang Tian\thanks{The second author's research is partially supported by an NSF grant and Simons fund. }\\
 Department of Mathematics\\
 Massachusetts Institute of Technology\\
 77 Massachusetts Avenue\\
 Cambridge, MA02139\\
 email:\texttt{tian@math.mit.edu}\\[3pt]}
\date{February 18, 2004}
\maketitle
\begin{abstract}In this paper, we prove a rigidity theorem of
asymptotically hyperbolic manifolds only under the assumptions on
curvature. Its proof is based on analyzing asymptotic structures
of such manifolds at infinity and a volume comparison theorem.
\end{abstract}
\section{Introduction}
In this paper, we study the rigidity problem for asymptotically
hyperbolic manifolds. Much progress has been made on this problem.
In \cite{M}, using the Dirac operator, Min-oo proved that a spin
manifold of dimension $n$ must be a hyperbolic space if it is
asymptotic to hyperbolic space in a strong sense and its scalar
curvature is not less that $-n(n-1)$. His argument was refined and
new exciting results were obtained by Andersson and Dahl \cite{AD}
and X.Wang \cite{W}. For even dimensional manifolds, Leung proved
in \cite{L} that any conformally compact Einstein manifold $({\Bbb
B}^n, g)$  which is asymptotically hyperbolic of order greater
than $2$ must be hyperbolic. By exploring properties of positive
eigenfunctions, J.Qing proved that a conformally compact Einstein
manifold with round sphere as its conformal infinity has to be a
hyperbolic space when the dimension is not greater than $7$ (cf.
\cite{Q}). He did not need to assume that the manifold considered
is spin. However, his approach relies on the positive mass theorem
for asymptotically flat manifolds. In all above results, one needs
to assume that there are nice coordinates at infinity and in such
coordinates, the metrics tensor behaves well. In view of geometry,
it would be natural to ask whether such an assumption can be
replaced by an intrinsic geometric condition. In this paper, we
will show a rigidity theorem of this type only under the
assumption on curvature.

Let $(X^{n+1},\ g)$ be a complete noncompact Riemannian manifold,
we call it an asymptotically locally hyperbolic manifold, which we
abbreviate as ALH in the following, of order $\alpha$ if
 $|K(x)+1|=O(e^{-\alpha \rho(x)})$, where $K(x)$ is the
sectional curvature of $g$ at point $x$ in any direction and
$\rho(x)=dist_g (x,o)$.

Recall that a Riemannian manifold $X$ has a pole $o$ if the
exponential map $exp_o:T_o X \rightarrow
 X$ is a diffeomorphism. Without loss of generality, in our case,
we may assume that the sectional curvature is negative outside a
unit ball of $(X,g)$. We have:

\begin{theo}
\label{th:main} Suppose that $(X^{n+1},g)$ $n\geq 2$ and $n\neq 3$
is an ALH manifold of order $\alpha$ with a pole and there is a
$\rho > 1$ such that the geodesic sphere with radius $\rho$ and
center at the pole is convex. If we further have $\alpha>2$ and
$Ric(g)\geq -ng$, then $(X^{n+1},g)$ is isometric to
$\Bbb{H}^{n+1}$.
\end{theo}

As a corollary, we have:
\begin{coro}
\label{th:maincorollary} Suppose that $(X^{n+1},g)$ $n\geq 2$ and
$n\neq 3$ is an ALH manifold of order $\alpha$ ($\alpha > 2$),
$K\leq 0$ and $Ric(g)\geq -ng$, then $(X^{n+1},g)$ is isometric to
$\Bbb{H}^{n+1}$.
\end{coro}

Let  $Rm^0$ denotes the traceless part of the curvature
tensor\footnote{The metric $g$ is of constant sectional curvature
iff $Rm^0$ vanishes. This property determines $Rm^0$
uniquely.},$\|Rm^0\|$ denote the norm of the tensor for $(X,g)$,
then for $n=3$, we have:
\begin{theo}
\label{th:main}Suppose that $(X^4,g)$ is an ALH manifold of order
$\alpha>2$ with a pole and there is a $\rho > 1$ such that the
geodesic sphere with radius $\rho$ and center at the pole is
convex. If we further have  $\|Rm^0 \| \in L^1 (X)$ and
$Ric(g)\geq -3g$,  then $(X^4,g)$ is isometric to $\Bbb{H}^4$.
\end{theo}

We will use the volume comparison theorem to prove above theorem.
In order to use the volume comparison, we need to estimate the
volume growth of geodesic spheres at infinity. We will carry this
out in several steps. First, we show that by changing the metric
conformally, we can compactify $(X,g)$ in an appropriate way.
Next, we will show that the boundary of compactified Riemannian
manifold is isometric to standard sphere, in this step, we first
verify that the boundary is conformal to the standard sphere. It
follows from the assumption on curvature that the boundary is
diffeomorphic to the standard sphere, hence, it suffices to show
that the boundary is locally conformally flat. By a direct
computation, we can show that the Weyl tensor of the boundary
vanishes, if the induced metric on the boundary is sufficiently
smooth we know that it is locally conformally flat. However, since
the metric on the compactified boundary is not necessarily smooth
enough, we have to check what the locally conformal flatness of
the boundary means in our current case. Under the assumption on
Ricci curvatures in above theorem, we observe that the scalar
curvature and volume of the boundary of compactified manifolds is
less or equal to those of the standard sphere. It follows that the
scalar curvature of the boundary is actually equal to that of the
standard sphere, hence, if $n=2$, we see that the boundary is
isometric to the standard sphere; if $n\geq 3$, then by Obata's
theorem, we know that the boundary is also isometric to the
standard sphere; Finally, we can show that the volume of geodesic
spheres of $(X, g)$ is equal to that of corresponding geodesic
spheres in $\Bbb{H}^{n+1}$ with the same radius, then, by the
volume comparison theorem, we prove the main theorem.

This assumption $\alpha>2$ should be optimal, since there are many
asymptotically hyperbolic Einstein metrics on $\Bbb{B}^4$ with
$\alpha=2$, we refer the readers to Theorem C and Appendix in
\cite{A} for details.  In the case of $n=3$, in order to show
locally conformal flatness of the boundary, one has to check that
certain linear combination of covariant derivatives of Schouten
tensor vanishes, for time being, we do not know how to deduce this
one from the assumption $\alpha >2$, this is the reason why we
need the extra assumption  $\|Rm^0\|_g \in L^1 (X)$, we doubt its
necessity. We also think that the assumption on existence of a
pole is unnecessary. In order to remove the assumption on pole,
one may study asymptotics of certain eigenfunctions at infinity
and use appropriate power of them to scale metrics as we do in the
next section. We will discuss this in a future paper. also one can
generalize arguments to study rigidity of asymptotic symmetric
spaces, one particularly interesting case is for asymptotic
complex hyperbolic K\"ahler manifolds. We expect that a similar
result can be proved for K\"ahler manifolds by assuming that
bisectional curvature tends to -1 at a sufficiently fast rate.

The organization of this paper is as follows: In Section 2, we
discuss the compactification and conformal structure of $(X,g)$ at
infinity; In Section 3, we show that the boundary of the
compactified manifold is isometric to the standard sphere and then
use it to deduce the main theorem.

{\bf Acknowledgment}: Part of this work was done during the first
author's visit at Department of Mathematics of MIT. He would like
to thank colleagues there for providing an excellent research
environment. Especially, he wants to thank Dr. X.Wang for
stimulating discussions during this visit.

\section{Compactfication and conformal structure at infinity}

In this section, we give a compactfication of $(X,g)$ at infinity
and study the induced conformal structure at infinity. This
compactification is crucial in the proof of our main theorem.

Let $\Sigma_\rho$ be the geodesic sphere in $(X,g)$ with radius
$\rho$ and a fixed center $o$. Define $\bar{g}$ to be
$\sinh^{-2}\rho g$, then we have:

\begin{theo}
\label{th:compactification} There is a subsequence of
$(\Sigma_\rho, \bar{g}_\rho)$ which converges to a $W^{2,p}\cap
C^{1,\alpha}$ Riemannian manifold $(\Sigma_\infty,
\bar{g}_\infty)$ in the weakly $W^{2,p}$-topology, where $p\in
(1,\infty)$ and $\alpha\in (0,1)$ are arbitrary. Here by a
$W^{2,p}\cap C^{1,\alpha}$ structure on $(\Sigma_\infty, \bar
g_\rho)$, we mean that there is a covering $\{U_i\}$ of
$\Sigma_\infty$ by coordinates $\phi_i: U_i\mapsto \RR^n$ such
that the transition functions $\phi_i\cdot \phi_j^{-1}$ and the
metric tensors $\phi_i^{-1*}g$ are in $W^{2,p}\cap C^{1,\alpha}$.
Furthermore, $(\Sigma_\infty, \bar{g}_\infty)$ is conformally
equivalent to the standard sphere. Here $\bar{g}_\rho$ denotes the
restriction of $\bar{g}$ to $\Sigma_\rho$.
\end{theo}

By the compactness theorem proved  in \cite{GW}, in order to have
the convergence property of $(\Sigma_\rho, \bar{g}_\rho)$, we only
need to show the following

\begin{lemm}
\lab{lemm:curvature-bound} There exists a constant $C$ such that
$|{Rm}(\bar g_\rho)|\leq C$, $Vol(\Sigma_\rho, {\bar g}_\rho) \geq
C^{-1}$, ${\rm diam} (\Sigma_\rho, {\bar g}_\rho) \leq C$, where
${Rm}(\bar g)$ denotes the curvature tensor of $\bar g$.
\end{lemm}

Let us first recall some basic formulae. For time being, we assume
$\bar g$$=$$u^2 g$. Let $\{\omega_a \}_{1\leq a \leq n+1}$ be a
local orthonormal coframe of $g$ such that $\omega_{n+1}=d\rho$
and $\{\omega_i \}_{1\leq i \leq n}$ is tangent to $\Sigma_\rho$.
For convenience, we also denote $g=d\rho^2 +g_{ij}(\rho,
\theta)d\theta_i d\theta_j$. Then we have structure equations

\begin{equation}
\left\{\begin{array}
     {r@{}l}
    d\omega_a & = \sum_{b=1}^{n+1} \omega_{ab}\wedge \omega_b,
    \quad
     \omega_{ab}+ \omega_{ba}=0 \\
    d\omega_{ab}&=\sum_{c=1}^{n+1}\omega_{ac}\wedge \omega_{cb} -
    \frac{1}{2}\sum_{c,d=1}^{n+1}R_{abcd}\omega_c \wedge \omega_d,
    \end{array}\right.
\end{equation}
where $R_{abcd}$ denote components of the curvature tensor. The
second fundamental form, denoted by $A=(h_{ij})_{1\leq i, j\leq
n}$, of $\Sigma_\rho$ with respect to $g$ is given by
$$\omega_{n+1 i}|_\Sigma =\sum_{j=1}^n h_{ij}\omega_j,
\quad h_{ij}=h_{ji},$$ where $(\cdot)|_\Sigma$ denotes the
restriction of an 1-form to $\Sigma_\rho$. The corresponding mean
curvature is given by $H=\sum_{i=1}^{n} h_{ii}$.

Let $\eta_a= u \omega_a $ ($1\leq a \leq n+1$), then $\{\eta_a
\}_{1\leq a\leq n+1}$ is an orthonormal coframe for the metric
$\bar g$, and

\begin{equation}
\left\{\begin{array}
     {r@{}l}
    d\eta_a & = \sum_{b=1}^{n+1} \eta_{ab}\wedge \eta_b,
    \quad
     \eta_{ab}+ \eta_{ba}=0 \\
    d\eta_{ab}&=\sum_{c=1}^{n+1}\eta_{ac}\wedge \eta_{cb} -
    \frac{1}{2}\sum_{c,d=1}^{n+1}{\bar R}_{abcd}\eta_c \wedge \eta_d
    \end{array}\right.
\end{equation}
where ${\bar R}_{abcd}$ are components of the curvature tensor of
$(X,\bar g)$ in the coframe $\{\eta_a \}_{1\leq a\leq n+1}$. By a
direct computation, we see that

$$\eta_{ab}=\omega_{ab}-(\log u)_b \omega_a + (\log u)_a \omega_b.$$
Here for any smooth function $f$ on $X$, $f_a$ is defined by
$df$$=$$\sum_{a=1}^{n+1}f_a \omega_a$. Thus, we get

$$\eta_{n+1, i}|_{\Sigma_\rho}=(h_{ij} + \frac{\partial}{\partial
\rho}(\log u)\delta_{ij})u^{-1}\eta_j.
$$
It follows that the second fundamental form of $\Sigma_\rho$ with
respect to $\bar g$ and $\{\eta_i\}_{1\leq i \leq n}$ is given by
\begin{equation}
\bar{h}_{ij}=(h_{ij} + \frac{\partial}{\partial \rho}(\log
u)\delta_{ij})u^{-1}.
\end{equation}

On the other hand, we can deduce from the structure equations for
curvatures
\begin{equation}\frac{\partial h_{ij}}{\partial \rho} +
\sum_{k=1}^n h_{ik}h_{kj} =- R_{n+1 i n+1 j}.
\end{equation}

In order to estimate $h_{ij}$, we need the following

\begin{lemm} Suppose that $f$ is a smooth function and for any $\rho >0$, we have
$|f(\rho)-1|\leq Ke^{-\alpha \rho}$ for some $\alpha>2$ and
$\frac14 \leq f(\rho) $. If $y$ is a solution of the equation
$$y' + y^2 = f(\rho)~~{\rm and}~~y(0)>0.
$$
Then there is a constant $C>0$, which depends only on $K$ and
$y(0)$, such that
$$|y-1|\leq C e^{-2\rho}. $$
\end{lemm}
\begin{proof} We will prove this lemma in the following steps.

\noindent {\bf Claim 1}: $0<y(\rho)\leq \rho +C_1$ for any
$\rho>0$.

Here and in the sequel, $C_i$ always denotes a constant which
depends only on $y(0)$ and $K$. Clearly, $y(\rho)\leq \rho +C_1$.
To see that $y(\rho)>0$, we first observe that $f\ge frac14$,
hence, by using the equation, $y'(\rho)> 0$ whenever $y(\rho) <
frac12$. It follows that $y$ increases in the region where
$y<\frac12$. Then the claim follows from $y(0)>0$.

\noindent {\bf Claim 2}: $|y-1|\le C_2 e^{-\rho}$.

Set $v=y-1$, we have $\rho +C_1 -1 \geq v \geq -1$ and $-1\le
v(0)\le y(0)$. Choose $\beta = 1+frac{\alpha}{2} < \alpha$. Then
$|v|\le C_3 e^{(\alpha -\beta )\rho}$, consequently, using the
equation for $y$, we can deduce
$$(v^2)' + 2 v^2 \le (v^2)' +(4+2v)v^2 \leq C_4 e^{-\beta \rho}, \quad 2<\beta<\alpha.$$
It follows
$$(v^2 e^{2\rho})' \leq C_4 e^{(2-\beta)\rho},$$
Integrating this inequality, we get
$$v^2 \le (v^2(0) + \frac{2 C_4}{\alpha - 2}) e^{-2\rho} \le C_5 e^{-2\rho}.$$
Claim 2 follows.

Now we can finish the proof of this lemma. By Claim 2, we have
$$|v|\leq C_5 e^{-\rho}.$$
Using this and the equation for $y$, we have
$$(v^2)' + (4 - 2 |v|) v^2 \le 2 K e^{-\alpha \rho} |v|.$$
Suppose that we have proved $v^2\le e^{-2\beta_k}$ for some
$\beta_k \ge 1$, then it follows from the above
$$(v^2 e^{4\rho})' \le C_6 (e^{(4-\alpha - \beta_k)\rho} + e^{(4-3\beta_k)\rho}.
$$
Integrating this, we get
$$v^2 \le C_7 (e^{-4\rho} + e^{-\min\{3\beta_k, \alpha + \beta_k\} \rho}).$$
If $\min\{3\beta_k, \alpha + \beta_k\}\ge 4$, we are done,
otherwise, then we take $\beta_{k+1}= \frac12 \min\{3\beta_k,
\alpha + \beta_k\} \ge \beta_k + \frac{\alpha}{2} - 1$ and repeat
the above process. Then the lemma follows after finitely many
iterations.
\end{proof}

\begin{lemm} Let $A$$=$$\sum_{ij}h_{ij}\omega_i \otimes \omega_j$ be the
second fundamental form of $\Sigma_\rho$ in $(X,g)$ and write
$$h_{ij}=\delta_{ij} +T_{ij}e^{-2\rho},$$
then $\|T\|_g \leq C$$<+\infty$, where $T=\sum_{ij}T_{ij}\omega_i
\otimes \omega_j$.
\end{lemm}
\begin{rema} If $\hat {h}_{ij}$ denotes components of the second
fundamental form of $\Sigma_\rho$ in $(X,g)$ in the coordinate
frame $\{\frac{\partial}{\partial \theta^i}\}$, then we have
$$\hat {h}_{ij} =g_{ij}+p_{ij} e^{-2\rho}.$$
Write $\omega_i =$$\sum_j b_{ij}d\theta^j$, we have
$(g_{ij})=$$(b_{ij})^T \cdot(b_{ij}) $ and $(p_{ij})=$$(b_{ij})^T
\cdot$ $(T_{ij}) \cdot(b_{ij})$.
\end{rema}

\begin{proof} Let $\lambda_{max}$ and $\lambda_{min}$ be the
largest and smallest eigenvalue of matrix $(h_{ij})$, then they
are Lipschitz, and we claim that
\begin{equation}
\lab{eq:equforeigenvalus} \frac{d}{d\rho}\lambda_{max} +
\lambda_{max}^2 =1+O(e^{-\alpha \rho}).
\end{equation}

\begin{equation} \frac{d}{d\rho}\lambda_{min} + \lambda_{min}^2
=1+O(e^{-\alpha \rho}).
\end{equation}

In fact, for any $\rho=\rho_0$, let $V$ be the unit eigenvector of
$\lambda_{max}$, then $V^T(h_{ij})V|_{\rho=\rho_0}$\quad
$=$$\lambda_{max}(\rho_0)$, and $V^T(h_{ij})V \leq
\lambda_{max}(\rho)$ for any $\rho$, thus,

$$\frac{d}{d\rho}\lambda_{max}|_{\rho=\rho_0}=\frac{d}{d\rho}V^T(h_{ij})V|_{\rho=\rho_0},$$
hence,

$$\frac{d}{d\rho}\lambda_{max}|_{\rho=\rho_0} + \lambda_{max}^2|_{\rho=\rho_0}
=1+O(e^{-\alpha \rho}),$$

which implies (2.2) is true, by the same reason, (2.3)is true too.
On the other hand, when $\rho$ is sufficiently large the
eigenvalue of matrix $(R_{n+1in+1j})$ is less than $-\frac14$, and
note that there is a convex geodesic sphere with sufficiently
large radius, hence, we may assume the initial data of equation
(2.3),(2.4) is positive, then by Lemma 2.3, we see Lemma 2.4 is
true.

\end{proof}

Due to Lemma 2.4, we have $\sup_{\Sigma_\rho}|T_{ij}|$$\leq C
<+\infty$ for any $\rho\geq 1$.

\noindent {\bf Proof of Lemma 2.2}: By a direct computation, we
have
\begin{equation}
\begin{split}\bar{R}_{abcd}&=u^{-2}R_{abcd}-u^{-2}((\log
u)_{bm}-(\log u)_m (\log
u)_b)(\delta_{ac}\delta_{dm}-\delta_{ad}\delta_{mc})\\
&-u^{-2}((\log u)_{am}-(\log u)_m (\log
u)_a)(\delta_{mc}\delta_{db}-\delta_{md}\delta_{bc})\\
&-u^{-2}|\nabla \log u|^2
(\delta_{ac}\delta_{bd}-\delta_{ad}\delta_{bc}),
\end{split}
\end{equation}
where $\bar{R}_{abcd}$ denote the components of the curvature
tensor of $(X,\bar g)$ in the coframe $\{\eta_a\}$. By our
assumption on asymptotic hyperbolicity, we may write
$$R_{abcd}=(\delta_{bc}\delta_{ad}-\delta_{ac}\delta_{bd})+E_{abcd},$$
where $|E_{abcd}|=O(e^{-\alpha \rho})$.

Now let $u=\sinh^{-1}\rho$. Noticing that for any $1\le a,b\le n$
$$(\log u)_{ab}=(\log u)_\rho (\delta_{ab}+
T_{ab} e^{-2\rho}),$$ we can deduce from the above and Lemma
2.4\footnote{Without loss of generality, we may assume that
$\alpha \le 4$.}

\begin{equation}
\bar{R}_{abcd}=\frac14(T_{bd}\delta_{ac}-T_{bc}\delta_{ad}
-T_{ad}\delta_{bc} + T_{ac}\delta_{bd})+E_{abcd}, \quad 1\leq
a,b,c,d \leq n
\end{equation}
here, $|E_{abcd}|=O(e^{(2-\alpha)\rho})$ as $\rho$ tends to
infinity.

$$\bar{R}_{n+1bcd}=O(e^{(2-\alpha)\rho}),$$

$$\bar{R}_{n+1bn+1d}=T_{bd}-\frac12\delta_{bd}
+O(e^{(2-\alpha)\rho}).$$ Let $\bar{h}_{ij}$ be the components of
the second fundamental form of $\Sigma_\rho\subset X$ with respect
to the metric $\bar g$. It follows from (2.3) that:

\begin{equation}
\bar{h}_{ij}=O(e^{-\rho}).
\end{equation}

Now let us estimate the volume and diameter of $(\Sigma_\rho, \bar
g)$. We can write $g$ in the form $d\rho^2 +g_{ij}(\rho,
\theta)d\theta^i d\theta^j$, then we have
$$\frac{\partial}{\partial \rho}g_{ij}=2{\hat h}_{ij}.$$
Using the facts that $\hat h_{ij} = g_{ij} + p_{ij} e^{-2\rho}$
and $-c (g_{ij}) \le (\hat h_{ij}) \le c (g_{ij})$ for some
constant $c$, we can show that there exists a constant $\Lambda$
independent of $\rho$ such that
\begin{equation}
\Lambda^{-1}e^{2\rho}(\delta_{ij}) \leq (g_{ij})\leq \Lambda
e^{2\rho}(\delta_{ij}).
\end{equation}
It follows that ${\rm diam}(\Sigma_\rho, \bar g)\leq C_2$ and
$Vol(\Sigma_\rho, \bar g)\geq \delta_0 >0$. The proof of Lemma 2.2
is completed.

By using (2.8), (2.9) and the Gauss equations, we see that the
sectional curvature of $(\Sigma_\rho, \bar g)$ is uniformly
bounded. Then it follows from Lemma 2.2 and \cite{GW} that there
exists a sequence of $(\Sigma_{\rho_i},\bar{g}_{\rho_i})$, which
will be denoted by $(\Sigma_{i},\bar{g}_{i})$, converges to
$(\Sigma_\infty, \bar{g}_\infty)$ in the sense of weak topology of
$W^{2,p}$ for any $p<\infty$, and for any $q$ $\in$
$\Sigma_\infty$, there is a coordinate charts $(B_q,\theta^i)$ in
which the components of $\bar{g}_\infty$ is $C^{1,\alpha}\bigcap
W^{2,p}$, $\forall p <+\infty$ and the curvature of
$(\Sigma_\infty, \bar{g}_\infty)$ is bounded.

Let $\hat{R}_{ijkl}$ be the components of the curvature tensor of
$(\Sigma_\rho, \bar g_\rho)$ under the orthonomal frame $\eta_i$
($1\leq i\leq n$), then the Weyl tensor is:
$$
\hat W_{ijkl}=\hat R_{ijkl}-\frac{1}{n-2}(\hat
R_{ik}\delta_{jl}-\hat R_{jk}\delta_{il} +\hat R_{jl}\delta_{ik}
-\hat R_{il}\delta_{jk})+ \frac{\hat
R}{(n-1)(n-2)}(\delta_{ik}\delta_{jl}-\delta_{jk}\delta_{il})
$$
Combined with (2.8) and (2.9), we see that $\|\hat W \|=o(1)$ as
$\rho$ tends to $\infty$, and the Ricci tensor is of the form
$$\hat R_{ij}=\frac14 [(n-2)T_{ij}+tr_{\bar g}T \delta_{ij}]+E_{ij},$$
where $E_{ij}=\bar g^{kl} E_{ikjl}$ and $|E_{ij}|=o(1)$ as $\rho$
tends to infinity.

Recall that the Schouten tensor of $g$ is
\begin{equation}\hat S_{ij}=\frac{1}{n-2}(2\hat R_{ij}-\frac{\hat
R}{n-1}\delta_{ij}).
\end{equation}
Therefore, the Weyl tensor of $(\Sigma_\infty, \bar g_\infty)$
vanishes in the $L^p$-sense. Together with Gauss equations and
Codazzi equations and (2.11), we  deduce
\begin{equation}
\begin{split}\bar \nabla_k \hat S_{ij}-\bar \nabla_j \hat
S_{ik}&=\frac12 e^{2\rho}\sinh \rho R_{n+1 kij}+\frac{2}{n-2}(\bar
\nabla_k E_{ij}-\bar \nabla_j E_{ik})\\
&-\frac{1}{(n-1)(n-2)}(\bar \nabla_k E \delta_{ij}-\bar \nabla_j E
\delta_{ik})
\end{split}
\end{equation}
If $n=3$, by the assumption that $\|Rm^0\|_g\in L^1 (X)$, we see
that there are $\rho_i$ which tend to infinity such that
$$\int_{\Sigma_{\rho_i}}||Rm^0||_g \sinh^3{\rho_i}  \rightarrow 0,$$
in particular, we have
$$\int_{\Sigma_{\rho_i}}|R_{4ijk}|e^{3{\rho_i}}  \rightarrow
0.$$ It follows that for any $\phi \in
C^{\infty}(\Sigma_{\rho_i})$,
\begin{equation}
\int_{\Sigma_{\rho_i}}\phi R_{4ijk}e^{2\rho_i}\sinh \rho_i
\rightarrow 0.
\end{equation}

Without loss of generality, we may assume $(\Sigma_{\rho_i}, \bar
g_{\rho_i})$ converges to $(\Sigma_\infty, \bar g_\infty)$, for
simplicity, in the sequel, $(\Sigma_\infty, \bar g_\infty)$ will
be denoted by $(\Sigma, \bar g)$ and the components of its
curvature tensor will be simply denoted by $\bar R_{ijkl}$. Then
if $n=3$, we see that the Schouton tensor of $(\Sigma, \bar g)$
satisfies the following equations in the sense of distribution,
\begin{equation}\bar \nabla_k \bar S_{ij}-\bar \nabla_j \bar
S_{ik}=0,
\end{equation}
that is, for any $\phi \in C^\infty(\Sigma)$, we have
$$\int_\Sigma \bar S_{ij}\bar \nabla_k \phi -\bar S_{ik}\bar \nabla_j \phi =0,$$
where $\bar \nabla_k$ are covariant derivatives of $(\Sigma, \bar
g)$ with respect to an orthonormal basis $\{e_i\}_{1\leq i \leq
n}$.

Now, we are in the position to show the following:
\begin{theo}
Suppose that $(\Sigma,g)$ is an $n$-dimensional Riemannian
manifold and the metric $g$ is in $W^{2,p}$ for any $1<p<\infty$.
If its curvature tensor is bounded and Weyl tensor $W=0$ if $n>3$;
and (2.14) is true if $n=3$, then $ g$ is locally conformally
flat, i.e., for any point $q\in V$, there is a neighborhood $U$,
such that in $U$, we have a positive $W^{2,p}$ function $f$ with
$g =f g_{\rm euc}$, where $g_{\rm euc}$ denotes a flat metric on
$U$.
\end{theo}

Clearly, it is a local result, hence, we need to consider only the
problem in a local coordinate chart, i.e., we assume that $\Sigma$
is a ball $B^n\subset \RR^n$ and $g=g_{ij}dx_idx_j$, where
$x_1,\cdots,x_n$ are euclidean coordinates of $\RR^n$. By our
assumption, $g_{ij}$ are $W^{2,p}$ functions on $B^n$ for any
$1<p<\infty$. It follows the the curvature tensor $R_{ijkl}$ and
the Christoffel symbol $\Gamma^i_{jk}$ are in $L^p$ and $W^{1,p}$
respectively. Hence, we can define covariant derivatives of
$R_{ijkl,h}$ in the sense of distribution, that is, for any $\phi
\in C^\infty _0 (B^n)$, we have:
$$\int_{B^n} R_{ijkl,h}\phi\sqrt {{\rm det}(g)} dx= -\int_{B^n} R_{ijkl}\frac{\partial}{\partial x^h}
(\phi \sqrt {\rm det}(g))dx +\int_{B^n}Rm*\Gamma \phi \sqrt {{\rm
det}(g)}dx,
$$
where ${\rm det}(g)=det(g_{ij})$ and $Rm*\Gamma$ refers to a
bilinear form of $R_{ijkl}$ and $\Gamma^i_{jk}$. Since $R_{ijkl}$
are in $L^p$ and $\frac{\partial}{\partial x^h}(\phi \sqrt {{\rm
det}(g)})$ is in $C^{\alpha}$ for some $\alpha>0$, the right hand
side of the above equation is well defined. Similarly, we can
define $R_{ijkl,hm}$ in the sense of distribution, that is, for
any $\phi \in C^\infty _0 (B^n)$, we have:
$$\begin{array}{l}
\int_{B^n}R_{ijkl,hm}\phi \sqrt{{{\rm det}(g)}} dx =
\int_{B^n}R_{ijkl}\frac{\partial^2}{\partial x^h \partial
x^m}(\phi \sqrt {{\rm det}(g)})dx -\int_{B^n}Rm*\Gamma
\frac{\partial}{\partial x^m}(\phi \sqrt {{\rm det}(g)})dx
\\-\int_{B^n}Rm*\frac{\partial}{\partial x^m}(\Gamma \phi \sqrt
{{\rm det}(g)})dx +\int_{B^n}Rm*\Gamma*\Gamma\phi\sqrt {{\rm
det}(g)} dx.
\end{array}
$$

Now we have:
\begin{lemm} Suppose that $g \in W^{2,p}$ for some $p>1$, then in
the distributional sense, we have the second Bianchi identity
\begin{equation}
R_{ijkl,h}+R_{ijlh,k}+R_{ijhk,l}=0
\end{equation}
and
\begin{equation}
R_{ik,lt}=R_{ik,tl}+Ric*Rm.
\end{equation}
That is, for any $\phi \in C^{\infty}_0 (B^n)$, we have
$$\int_{B^n} (R_{ijkl,h}+R_{ijlh,k}+R_{ijhk,l})\phi \sqrt{{\rm det}(g)} dx =0$$
and
$$\int_{B^n}(R_{ik,lt}-R_{ik,tl}-Ric*Rm)\phi \sqrt {{{\rm det}(g)}}dx=0.$$
Here $R_{ij}$ is the Ricci tensor of $g$ and $Ric*Rm$ denotes a
bilinear form of Ricci tensor and curvature tensor.
\end{lemm}

\begin{proof}
By the assumption, we may take a sequence of smooth metrics $g_i$
which converges to $g$ in $W^{2,p}$. Since (2.15) and (2.16) hold
for the curvature tensor of $g_i$ and curvature tensors of $g_i$
converge to that of $g$ in $L^p$, we see that (2.15) and (2.16)
hold for $g$, too.
\end{proof}

Next we construct harmonic coordinates around any point of
manifold. Without loss of generality, we only need to show
\begin{lemm} Suppose that $g_{ij}$ are in $W^{2,p}$ on $B^n$ for any $1<p<\infty$,
then there are harmonic coordinates $(z^1, \cdots, z^n)$ around $o
\in B^n$ with $z^i$ in $W^{3,p}$.
\end{lemm}

\begin{proof} Let $\Gamma^i _{jk}$ denote the Christoffel symbols
of $g$ in euclidean coordinates $x^1,\cdots, x^n$. Define $y_i$ by
$x^i =y^i -\Gamma^i _{jk}(o)y^j y^k$ ($1\leq i \leq n$), by the
Inverse Theorem, we see that $y^i$ are smooth functions of
$(x^1,\cdots, x^n)$ around $o$ and form coordinates. Let $\bar
\Gamma^i _{jk}$ be the Christoffel symbols of $g$ in coordinates
$(y^1, \cdots, y^n)$, then by direct computations, we see
$$\bar \Gamma^k _{ij}\frac {\partial x^s}{\partial y^k} =
\frac{\partial ^2 x^s}{\partial y^i \partial y^j}+ \frac{\partial
x^l}{\partial y^j}\frac{\partial x^m}{\partial y^i}\Gamma^s_{lm}.
$$
It follows that $\bar \Gamma^k_{ij}(o)=0$, consequently, $\Delta
y^i=0$ at $o$. This implies that $||\Delta y^i||_{L^\infty
(B_\epsilon (o))}$ tends to $0$ as $\epsilon$ goes to zero.
Consider the following boundary value problem on $B_\epsilon (o)$
\[ \left\{\begin{array}
     {r@{}l}
    \Delta z^i& =0  \\
    z^i |_{\partial B_\epsilon}&=y^i |_{\partial B_\epsilon}.
  \end{array}\right.
\]
then, using standard estimates for elliptic equations, we get
$$||z^i-y^i||_{C^{1,\alpha} }\leq C ||\Delta y^i||_{L^\infty (B_\epsilon)}.$$
Therefore, $z^1, \cdots, z^n$ form local coordinates on
$B_\epsilon (0)$ when $\epsilon$ is sufficiently small. Clearly,
$z^i$ are in $W^{3,p}(B_\epsilon)$ and harmonic with respect to
$g$. The lemma is proved.
\end{proof}

By a direct computation and Lemma 2.8, we see that metric tensor
of $g$ in coordinates $z^1,\cdots,z^n$ is also in $W^{2,p}$. In
the following, we will consider the problem in these harmonic
coordinates, and the metric components will be still denoted by $
g_{ij}$.

\begin{lemm}
Let $R$ be the scalar curvature of $g$ and bounded, then when
$\epsilon$ is sufficiently small, the following equation
\[ \left\{\begin{array}
     {r@{}l}
    \Delta u-\frac{n-2}{4(n-1)}R u& =0  \\
    u |_{\partial B_\epsilon}&= 1|_{\partial B_\epsilon}.
  \end{array}\right.
\]
has a positive solution in $W^{2,p}(B_\epsilon)$.
\end{lemm}

\begin{proof} We note that when $\epsilon$ is sufficiently small,
the first Dirichlet eigenvalue can be arbitrarily large, and $ R$
is bounded, hence, the corresponding homogenous equation has only
trivial solution, and  this implies the above equation has
nonnegative solution, then by Lemma 3.4 in \cite{GT} (p34), we see
that the solution has to be positive. This finishes the proof of
the lemma.
\end{proof}

In order to show Theorem 2.6, we need the following lemma ( see
Theorem 17.2.7, \cite{H}, p18 for its proof).

\begin{lemm} Let $a_{ij}(x)$ be Lipschitz continuous in an open set
$\Omega \subset {\Bbb R}^n$, and assume that the matrix $(a_{ij})$
is positive definite and $u\in L^2_{loc}(\Omega)$. Then
$$\Sigma \frac{\partial}{\partial x^j}(a_{jk}\frac{\partial u}{\partial x^k})=f,$$
implies $u\in W^{1,2}_{loc} (\Omega)$ if $f\in
H^{-1}_{loc}(\Omega)$, moreover, if $f\in L^2_{loc}(\Omega)$, then
$u\in W^{2,2}_{loc} (\Omega)$. Here $H^{-1}_{loc}(\Omega)$ is the
dual space of $W^{1,2}_{0} (\Omega)$.
\end{lemm}

Now we can finish the proof of Theorem 2.6. Since the scalar
curvature of $(\Sigma,g)$ is bounded, by Lemma 2.9, we may choose
a sufficiently small neighborhood of $q$ such that there is a
positive $W^{2,p}$ function $u$ on this neighborhood such that the
scalar curvature of $\bar g= u^{\frac{4}{n-2}}g$ vanishes. It is
easy to show that $\bar g$ is also in $W^{2,p}$ for any
$1<p<\infty$, moreover, its Weyl tensor also vanishes if $n\geq 4$
and (2.14) still holds if $n=3$. By Lemma 2.8, we can choose
harmonic coordinates of the metric $\bar g$ with metric tensor
$\bar g_{ij}$ in $W^{2,p}$. It suffices to show that the
corresponding Ricci tensor is smooth in these harmonic
coordinates. In the sequel, we will do everything in these
coordinates.

Since the Weyl tensor and the scalar curvature vanish, we have
\begin{equation}
{\bar R}_{ijkl}=\frac{1}{n-2}({\bar R}_{ik}{\bar g}_{jl}-{\bar
R}_{jk}{\bar g}_{il}+{\bar R}_{jl}{\bar g}_{ik} -{\bar
R}_{il}{\bar g}_{ik})
\end{equation}

On the other hand, by Lemma 2.7, we have the second Bianchi
identity for $\bar g$, hence, by a direct computation, we deduce
\begin{equation}
{\bar g}^{jh}{\bar R}_{jk,h}=0.
\end{equation}
If $n\geq 4$, using (2.17), (2.18) and the Bianchi identity, we
can also derive
\begin{equation}
{\bar R}_{il,k}-{\bar R}_{ik,l}=0.
\end{equation}
When $n=3$, since the scalar curvature vanishes, the above
equation is nothing but (2.14). It follows
$${\bar g}^{kt}{\bar R}_{il,kt}-{\bar g}^{kt}{\bar R}_{ik,lt}=0,$$
and because of (2.16) in Lemma 2.7, we have
$${\bar g}^{kt}{\bar R}_{ik,lt}={\bar g}^{kt} {\bar R}_{ik,tl}+{\bar g}*{\bar Ric}*{\bar Rm}.$$
Note that
$${\bar g}^{kt}{\bar R}_{ik,tl}=({\bar g}^{kt}{\bar R}_{ik,t})_{,l}=0,$$
so we have
$${\bar g}^{kt}{\bar R}_{il,kt}={\bar g}*{\bar Ric}*{\bar Rm}.$$
Since $\bar g$ is in $W^{2,p}$, the above equation can be written
as
\begin{equation}
\frac{\partial}{\partial x^t}({\bar g}^{kt} \frac{\partial{\bar
R}_{il}}{\partial x^k})={\partial {\bar g}}*{\partial {\bar
Ric}}+{\bar g}*{\bar Ric}*{\bar Rm}.
\end{equation}
Noticing that $\bar g \in W^{2,p}$ and $\bar Rm \in L^p$ for any
$p>1$, we see that the right hand side of (2.20) is in $H^{-1}$,
which is dual to $W^{1,2}_0$. Then it follows from Lemma 2.10 that
$\bar R_{ij}$ are actually in $W^{1,2}_{loc}$, in turns, this
implies that the right side of (2.20) is in $L^2_{loc}$, then
again by Lemma 2.10, we see that $\bar R_{ij}$ are in
$W^{2,2}_{loc}$, then it follows from the standard theory for
elliptic equations that $\bar R_{ij}$ are actually
$C^{2,\alpha}_{loc}$, therefore, $\bar g$ is smooth, and
consequently, by the classsical Weyl Theorem, it is locally
conformal flat. Theorem 2.6 is proved.

Now, we can prove Theorem 2.1.

\noindent {\bf Proof of Theorem 2.1}: It only remains to show that
$(\Sigma, \bar{g})$ is conformally equivalent to the standard
sphere. By the assumption of Theorem 2.1, we see that $\Sigma$ is
diffeomorphic to $\Bbb {S}^n$. On the other hand, by Theorem 2.6,
we know that $(\Sigma, \bar{g})$ is a locally conformally flat
manifold, so is conformally equivalent to $\Bbb{S}^n$. Theorem 2.1
is proved.

\section{Proof of Main  Theorems }

To prove Theorem 1.1 and Theorem 1.3, we need to compare both the
volume and the scalar curvature of $(\Sigma, \bar g)$ which is the
boundary of Riemannian manifold $(X, \bar g)$ with those
corresponding quantities of the standard sphere. Using this, we
are able to show that $(\Sigma, \bar g)$ is actually isometric to
the standard sphere. Then by the Volume Comparison theorem, we can
conclude that the original manifold $(X,g)$ is isometric to
$\Bbb{H}^{n+1}$.

\begin{lemm}
\lab{le:volume-scalar} Let $\omega_n$ denote the volume of
$\Bbb{S}^n$ and $\bar R$ be the scarlar curvature of $(\Sigma,
\bar{g})$, then we have ${\rm Vol}(\Sigma, \bar{g})\leq \omega_n$
and $\bar{R} \leq n(n-1)$.
\end{lemm}
\begin{proof} Recall that  $\hat {R}$ is the scalar curvature of
$(\Sigma_\rho, \bar{g}_\rho)$, then by the computations in last
section, we have
$$\hat{R}=\frac{n-1}{2}\sum_{i,j=1}^n g^{ij} p_{ij} + o(1),
\ \text{as}\ \rho\rightarrow \infty$$ and
$$H=n+e^{-2\rho} \sum_{i,j=1}^n g^{ij} p_{ij},$$
where $H$ denotes the mean curvature of $\Sigma_\rho$ in $(X,g)$.
On the other hand, because of $Ric(g)\geq -ng$, we can use the
Laplacian Comparison Theorem to get
\begin{equation}
H|_{\Sigma_\rho}=\triangle_g \rho |_{\Sigma_\rho} \leq n \coth
\rho
\end{equation}
It follows
$$\sum_{i,j=1}^n g^{ij} p_{ij} \leq \frac{2n e^{2\rho}}{e^{2\rho}-1},$$
and consequently,

$$\hat{R} \leq n(n-1) +o(1),$$
letting $\rho$ go to $\infty$, we get $\bar R \leq n(n-1)$.

To show ${\rm Vol}(\Sigma, \bar{g})\leq \omega_n$, we only need to
prove for any $\rho>0$,
\begin{equation}
{\rm Vol}(\Sigma_\rho, g)\leq (\sinh \rho)^n\, \omega_n.
\end{equation}
For any $\delta>0$, integrating (3.1), we obtain

$$\int_{B_{\tau +\delta}\setminus B_\tau}\triangle_g \rho dV_g \leq
\int_{B_{\tau +\delta}\setminus B_\tau} n \coth\rho dV_g,$$ which
is equivalent to

$$\frac{{\rm Vol}(\Sigma_{\tau+\delta}-{\rm Vol}(\Sigma_\tau)}{\delta} \leq
\frac n\delta\int_\tau^{\tau+\delta} \coth \rho {\rm
Vol}(\Sigma_\rho) d\rho.$$ Let $\delta\rightarrow 0$, we have

$$(\log \left((\sinh\tau)^{-n}{\rm Vol}(\Sigma_\tau)\right))'\leq 0.$$
Hence $(\sinh\tau)^{-n}{\rm Vol}(\Sigma_\tau)$ is non-increasing
with $\tau$. Since $\lim_{\tau\rightarrow 0}(\sinh\tau)^{-n} {\rm
Vol}(\Sigma_\rho) =\omega_n$, we see from the above that (3.2) is
true. This implies that ${\rm Vol}(\Sigma, \bar{g})\leq \omega_n$.
Thus Lemma 3.1 is proved.
\end{proof}

Our next goal is to establish

\begin{lemm} The limit space $(\Sigma, \bar{g})$ is isometric to
the standard sphere $(\Bbb{S}^n, g_0)$
\end{lemm}
\begin{proof} If $n=2$, we only need to show $\bar R=2$, suppose
not, we have $\bar R<2$ and $\rm {Vol}(\Sigma, \bar {g}) \leq
4\pi$, this is in contradiction with Gauss-Bonnet formula.

If $n\geq 3$, it suffices to prove that $\bar R =n(n-1)$. In fact,
we can write
\begin{equation}
\bar{g} =u^{\frac{4}{n-2}}g_0
\end{equation}
for some $u>0$ which belongs to $W^{2,p}$ for any $p<\infty$. If
$\bar R = n(n-1)$, $u$ satisfies a semi-linear elliptic equation
and the standard regularity theory implies that $u$ is smooth.
Then Lemma 3.2 follows from the Obata theorem.

Let $d \bar V$ and $dV_0$ be the volume elements of $(\Sigma,
\bar{g})$ and $(\Bbb{S}^n, g_0)$, respectively, then by (3.3),
$d\bar V= u^{\frac{2n}{n-2}}dV_0$. The following equation is
well-known

$$\bar R = u^{\frac{n+2}{2-n}}(n(n-1)u -\frac{4(n-1)}{n-2} \triangle_{\Bbb{S}^n} u),$$
it follows

$$\int_{\Bbb{S}^n} \bar{R} u^{\frac{2n}{n-2}}dV_0=
\int_{\Bbb{S}^n} ((n-1)nu^2+
\frac{4(n-1)}{n-2}|\nabla_{\Bbb{S}^n}u|^2)dV_0,$$ since $\bar{R}
\leq n(n-1)$, we get

$$n(n-1)(\int_{\Bbb{S}^n}u^{\frac{2n}{n-2}}dV_0)^{\frac2n}\geq
\frac{\int_{\Bbb{S}^n}((n-1)nu^2
+\frac{4(n-1)}{n-2}|\nabla_{\Bbb{S}^n}u|^2)dV_0}{(\int_{\Bbb{S}^n}u^{\frac{2n}{n-2}}dV_0)^{\frac{n-2}{n}}}.$$

Using the fact that $d \bar V =u^{\frac{2n}{n-2}}dV_0$ and
$Vol(\Sigma,\bar g) \leq \omega_n$, we see that:

\begin{equation}
n(n-1){\omega_n}^\frac2n \geq \frac{\int_{\Bbb{S}^n}((n-1)nu^2
+\frac{4(n-1)}{n-2}|\nabla_{\Bbb{S}^n}u|^2)dV_0}{(\int_{\Bbb{S}^n}u^{\frac{2n}{n-2}}dV_0)^{\frac{n-2}{n}}}.
\end{equation}

Write $g_0 = \psi^\frac{4}{n-2}ds_{{\Bbb R}^n}^2$, where
$\psi(x)=(\frac{1+|x|^2}{2})^{\frac{2-n}{2}}$, then $dV_0 =
\psi^\frac{2n}{n-2}dx$, where $dx$ is the volume element of
$\Bbb{R}^n$, we have:

\begin{equation}\text{The RHS of (3.4)}\ =\frac{4(n-1)}{n-2}\frac{\int_{{\Bbb
R}^n} |\nabla_{{\Bbb R}^n}(u \psi)|^2 dx}{(\int_{{\Bbb
R}^n}(u\psi)^\frac{2n}{n-2}dx)^{\frac{n-2}{n}}}.
\end{equation}

On the other hand, by a direct computation, we have:

\begin{equation}
n(n-1){\omega_n}^{\frac2n}=\frac{4(n-1)}{n-2}\frac{\int_{{\Bbb
R}^n} |\nabla_{{\Bbb R}^n}\psi|^2 dx}{(\int_{{\Bbb
R}^n}\psi^\frac{2n}{n-2}dx)^{\frac{n-2}{n}}}
\end{equation}

Putting (3.4), (3.5) and (3.6) together, we obtain
\begin{equation}
\frac{\int_{{\Bbb R}^n} |\nabla_{{\Bbb R}^n}\psi|^2
dx}{(\int_{{\Bbb R}^n}\psi^\frac{2n}{n-2}dx)^{\frac{n-2}{n}}} \geq
\frac{\int_{{\Bbb R}^n} |\nabla_{{\Bbb R}^n}u \psi|^2
dx}{(\int_{{\Bbb R}^n}(u\psi)^\frac{2n}{n-2}dx)^{\frac{n-2}{n}}}.
\end{equation}
Note that $\psi=(\frac{1+|x|^2}{2})^{\frac{2-n}{2}}$, we know that
the LHS of (3.7) is the best Sobolev constant for $\RR^n$, hence,
the equality in (3.7) holds, so $\bar R =n(n-1)$. Thus we see that
$(\Sigma, \bar{g})$ is nothing but $({\Bbb S}^n , g_0)$. Lemma 3.1
is proved.
\end{proof}

\noindent {\bf Proof of Theorem 1.1}: In the proof of Lemma 3.1,
we have shown that $(\sinh\rho)^n{\rm Vol}(\Sigma_\rho,g)$ is
non-increasing. By Lemma 3.2 and the fact that $(\Sigma_\rho,{\bar
g}_\rho)$ subconverges to $(\Sigma, \bar{g})$ in the
Cheeger-Gromov topology, we get
$$\lim_{\rho\rightarrow \infty}(\sinh\rho)^{-n}{\rm Vol}(\Sigma_\rho,g)
=\omega_n.$$ Hence, by the Volume Comparison Theorem, we have that
for any $\rho>0$, $(\sinh\rho)^{-n}{\rm
Vol}(\Sigma_\rho,g)=\omega_n$. Now we claim that
$$\triangle_g \rho =H|_{\Sigma_\rho} = n\coth \rho, \ \forall \rho>0.$$
If it is false, there is a point $p\in \Sigma_\rho$ such that
$\triangle_g \rho |_p < n\coth \rho|_p$, so
$$\int_{B_{\rho+\delta}(o)\setminus B_\rho (o)}\triangle_g \tau dV_g <
\int_{B_{\rho+\delta}(o)\setminus B_\rho (o)}n\coth \tau dV_g,$$
or equivalently
$${\rm Vol}(\Sigma_{\rho+\delta})-{\rm Vol}(\Sigma_\rho)<
n \int_{\rho}^{\rho+\delta}\coth \tau Area(\Sigma_\tau)d\tau.$$
This contradicts to that ${\rm Vol}(\Sigma_\tau)=(\sinh \tau
)^n\omega_n$. Hence for any $\rho>0$, $H|_{\Sigma_\rho}=n\coth
\rho$ and consequently
$$\frac{\partial H}{\partial \rho}+\frac{H^2}{n}=n.$$
However, from (2.2), we see that
$$\frac{\partial H}{\partial \rho}+\frac{H^2}{n}\leq n,$$
moreover, the equality holds if and only if $h_{ij}=\coth\rho
g_{ij}$. On the other hand, a direct computation shows that
$$\frac{\partial g_{ij}}{\partial \rho}=2 h_{ij}=2 \coth \rho g_{ij},$$
and
$$\lim_{\rho\rightarrow 0}\rho^{-2} g_{ij} = (g_0)_{ij},$$
Hence, $g_{ij}=(\sinh \rho )^2(g_0)_{ij},$ where $g_0$ is the
standard metric on $\Bbb{S}^n$. Therefore, we see $g= d\rho^2
+(\sinh\rho)^2 (g_0)_{ij}d\theta^i d\theta^j$, that is, $(X,g)$ is
isometric to ${\Bbb H}^{n+1}$. Theorem 1.1 is proved.

\end{document}